\documentclass[a4paper,twoside,english]{amsart}
\usepackage[T1]{fontenc}
\usepackage[latin9]{inputenc}
\usepackage{babel}
\usepackage{verbatim}
\usepackage{amsthm}
\usepackage{amsmath}
\usepackage{amssymb}
\usepackage[all]{xy}
\usepackage[unicode=true,
 bookmarks=true,bookmarksnumbered=true,bookmarksopen=false,
 breaklinks=false,pdfborder={0 0 1},backref=false,colorlinks=false]
 {hyperref}
\hypersetup{
 pdfauthor={P. Corvaja, A. Turchet and U. Zannier}}

\usepackage{tikz-cd}

\makeatletter

\pdfpageheight\paperheight 
\pdfpagewidth\paperwidth

\theoremstyle{plain}
\newtheorem{thm}{\protect\theoremname}[section]
 \newcommand\thmsname{\protect\theoremname}
 \newcommand\nm@thmtype{theorem}
 \theoremstyle{plain}

  \theoremstyle{remark}
  \newtheorem{rem}[thm]{\protect\remarkname}
  \theoremstyle{definition}
  \newtheorem*{example*}{\protect\examplename}
  \theoremstyle{definition}
  \newtheorem{example}[thm]{\protect\examplename}
  \theoremstyle{plain}
  
  \theoremstyle{plain}
  \newtheorem{prop}[thm]{\protect\propositionname}
  \theoremstyle{plain}
  \newtheorem{cor}[thm]{\protect\corollaryname}
  \newtheorem*{thm*}{\protect\theoremname}
 
\makeatletter
  \theoremstyle{definition}
  \newtheorem{my@rem}[thm]{Remark}
  \renewenvironment{rem}{\begin{my@rem}}{\end{my@rem}}
\makeatother

\makeatother

  \providecommand{\examplename}{Example}
  \providecommand{\lemmaname}{Lemma}
  \providecommand{\propositionname}{Proposition}
  \providecommand{\remarkname}{Remark}
  \providecommand{\theoremname}{Theorem}
\providecommand{\theoremname}{Theorem}
 \providecommand{\corollaryname}{Corollary}
\usepackage[left=3.2cm,right=3.5cm,top=3cm,bottom=3cm]{geometry}

\usepackage{mathtools}

\def\N{{\mathbb N}}
\def\Q{{\mathbb Q}}

\def\Z{{\mathbb Z}}
\def\G{{\mathbb G}}
\def\A{{\mathbb A}}

\def\O{{\mathcal O}}
\def\J{{\mathcal J}}

\def\S{{\mathcal S}}
\def\calC{{\mathcal C}}
\def\calF{{\mathcal F}}
\def\I{{\mathcal I}}
\def\Gal{\mathop{\rm Gal}\nolimits}

\def\P{{\mathbb P}}
\def\C{{\mathbb C}}

\def\CVD{{\hfill\hfil{\lower 2pt\hbox{\vrule\vbox to 7pt
{\hrule width  5pt\varphifill\hrule}\varphirule}}}\par}

\begin{document}

\author{Pietro Corvaja}
\address{Pietro Corvaja - Dipartimento di Matematica e Informatica\\
Universit\`a di Udine\\
Via delle Scienze, 206\\
33100 Udine --  Italy}
\email{pietro.corvaja@dimi.uniud.it}
\author{Amos Turchet}
\address{Amos Turchet - Dipartimento di Matematica e Fisica \\
Universit\`a di Roma Tre\\
Largo San Murialdo, 1\\
00146 Roma -- Italy}
\email{amos.turchet@uniroma3.it}

\author{Umberto Zannier}
\address{Umberto Zannier - Scuola Normale Superiore\\
Piazza dei Cavalieri, 7\\
56126 Pisa -- Italy}
\email{umberto.zannier@sns.it}
%---
%---
\title{Around the Chevalley-Weil Theorem}
\date{\today}

\begin{abstract}
  We present a proof of the Chevalley-Weil Theorem that is somewhat different from the proofs appearing in the literature and with somewhat weaker hypotheses, of purely topological type. We also provide a discussion of the assumptions, and an application to solutions of generalized Fermat equations, where our statement allows to simplify the original argument of Darmon and Granville.
\end{abstract}

\maketitle

\bigskip

\section{Introduction}

This short note does not claim substantially new results. Rather, its  main goal   is to illustrate  a self-contained proof of the Chevalley-Weil Theorem, originally stated in \cite{CW}, with a rather different presentation (compared to the usual ones) and  mild assumptions,  of purely topological content. To our knowledge, this does not appear in this form in the literature.

Hopefully this would  provide an accessible statement of the theorem,  ready-made for applications, without necessity of checking a number of requirements which sometimes are not quite the same in the various existing proofs and treatments. 

We shall also discuss and compare various concepts of {\it ramification}, obtaining as a byproduct some conclusions which appear to be new. 

\medskip

In order to describe in more detail the paper, let us start by recalling some generalities. 
The Chevalley-Weil Theorem is of course a precious tool in Diophantine Geometry: roughly speaking, given an {\it unramified}  morphism of finite degree  $\pi:W\to V$ of projective algebraic varieties defined over a number field $k$,  it allows 
 to {\it lift} the rational points $p\in V(k)$ to points of $W$ defined over a {\it fixed} number field, whereas in general if $\deg \pi>1$ and $V$ has enough rational points, we expect the field of definitions of the fibers $\pi^{-1}(p)$ to vary with $p$ so to generate a field of infinite degree. There is also a general version   for {\it integral } points on a quasi-projective variety. 
 
 The bulk of the assumption is that the morphism  should be {\it unramified}.  The result may then be seen as an arithmetical analogue of the familiar  {\it lifting theorem}  for maps,  in basic  homotopy theory. 

\medskip

The  known proofs  of this   lifting statement  split into two steps. The first one, which is the kernel of the theorem,  asserts that the extensions of residue fields of rational points are unramified (outside a suitable set of primes). The second  step appeals to a theorem of Hermite, predicting finiteness of the set of number fields with bounded degree and discriminant.  

This second part is entirely  independent of the first one, and has auxiliary nature. Therefore let us focus on the first step (which  is in fact the only part stated in the original paper). 

\medskip

There are several proofs available in the literature, and, beyond the original paper \cite{CW} (devoted only to curves),  we refer for instance to the books by Bombieri-Gubler \cite{BG}, by Lang \cite{L} and by Serre \cite{SeMWT}.  Most proofs use local discriminants, in order to exploit the assumption of absence of ramification (deducing it for the number fields obtained after specialisation). This (fundamental) assumption is usually formulated in the appropriate algebraic sense.

\medskip

In detail, our goal  is twofold:  on the one hand, we adopt a different (though  equivalent) viewpoint in exploiting non-ramification compared to other proofs. On the other hand, we adopt a purely topological notion of ramification. This will also  lead us to an apparently novel result, perhaps not entirely free of independent interest.

More precisely: 

\medskip

(i) An unramified cover  (of finite degree) may be viewed (by going to a normal closure) as a quotient by a group action without fixed points: we exploit this viewpoint rather than discriminants. The absence of fixed points reflects in the triviality of inertia groups of the specialised number fields, which amounts to non-ramification. In our presentation we adopt precisely this viewpoint and properties.  (Ultimately of course we use the same facts as in other proofs, but looking at  fixed points seems to us a somewhat different perspective.)  

\medskip

(ii)  There are several notions of {\it unramified morphism}. Algebraically, some definitions go back to the paper \cite{LS} of Lang and Serre, and   there is Grothendieck's notion of { \it \'etale morphism}.  On the other hand  of course there is the topological notion of covering space. There are also some implications and comparisons among these notions; but usually there is equivalence among the notions only under some supplementary assumption, as for instance {\it normality} of the relevant varieties.  In   Section \ref{sec:topo} we discuss the precise assumptions that are needed, and also prove a result showing that the absence of ramification for varieties defined over a number field does not depend on the embedding of the number field in $\C$. This is easy and well-known for normal varieties but seems not to have been stated before in general.

We hope that this discussion, though probably not needed for most applications, can be not free of interest for some readers, and may be even helpful for some applications, so to avoid the need for changing model of some variety or for checking certain assumptions. 

\medskip

We shall conclude this note with a direct application of Theorem \ref{thm:CW} to the study of solutions of generalized Fermat equations. In Section \ref{sec:appl} we provide a simplification of the original argument of Darmon and Granville in \cite[Theorem 2]{DG} where our Theorem \ref{thm:CW_integral}, applied to a quasi-projective variety which is not affine, allows to avoid the use of a  result of Beckmann \cite{Beck}.

\medskip

Let us now go to some precise statements. We recall that a continuous surjective map between topological spaces $\pi: W\to V$ is said to be a {\it topological cover} if the following holds: for every point $p\in V$ there exists an open neighborhood $U\subset V$ of $p$, a discrete  set $F$ and a homeomorphism $\pi^{-1}(U) \to U\times F$ making the following diagram to commute
\begin{equation*}
  \xymatrix{\pi^{-1}(U) \ar[rr] \ar[dr] & & U \times F \ar[dl] \\ & U & }
\end{equation*}
(where the arrow  $U\times F \to U$ denotes the projection on $U$).

In all our statements, algebraic varieties are allowed to possibly be reducible.

\begin{thm}[Chevalley-Weil Theorem]
  \label{thm:CW}
  Let $\pi:W\to V$ be a   dominant morphism of  %finite degree of   
  projective varieties over a number field $k$. Assume  there exists an embedding $\varphi: k \to \C$ such that the corresponding map $\pi_{\varphi,\C}=\pi_{\C}: W(\C) \to V(\C)$ is a  topological cover (i.e., viewing $V,W$ as varieties defined over $\C$ through $\varphi$). Then there exists a finite set $S$ of places of $\Q$ such that if $P\in V(\bar k)$, $Q\in W(\bar k)$, $\pi(Q)=P$, then the extension $k(Q)/k(P)$ is  ramified at most above places  of  $k(P)$ lying above  $S$.
\end{thm}

Note the purely topological formulation of the main hypothesis  
 (i.e., $\pi_\C$ is a topological covering map). 
In particular, we do not assume that the map $\pi$ is \'etale nor even that it is unramified in the algebraic sense (i.e. $\Omega^1_{W/V} = 0$) as in some  standard references for this result, e.g. \cite[Theorem 10.3.5]{BG}, \cite[Theorem 8.1]{L} or \cite[Section 4.2]{SeMWT}. In fact, our hypotheses apply in situations where the map $\pi$ is \emph{not} (algebraically) unramified (e.g. as in the normalization of the cuspidal cubic curve).
See Section \ref{sec:topo} for more examples and details.

In that section we shall also prove that the condition that $\pi_\C$ is a topological cover does not depend on the embedding $\varphi$ (see Proposition \ref{prop:conj}). This is well-known in most situation, that is, under quite general hypotheses, for instance when $V$ is normal. However we haven't  been able to locate a proof in the literature in the most general case.

\medskip

Several statements in the literature (but not the original one) assume that the point $P$ is defined over a fixed number field (and then, as mentioned above,  the above conclusion is combined with a theorem of Hermite to deduce that all points in the fiber over $P$ are defined over a(nother) fixed  number field). Here, instead, the field of definition of $P$ is arbitrary.

 \medskip

{\tt $S$-integral points, good reduction}. In the same setting there is also  an analogue result for integral points. We recall a few definitions for completeness. 

\medskip

For a finite set $S$ of places of $\Q$ containing the real place, we define $\O_{S}$ to be the subring of the algebraic closure $\overline\Q$ formed of elements which are integral at any place not lying above a place in $S$.  For a number field $k$, $\O_{k,S}$ is set to be $k\cap \O_S$.

For a quasi-projective variety $V$, we assume that an embedding in $\P_n$ is given for some $n$, for which a reduction of $V$ $\pmod v$ is defined for every place $v$ of $k$ not lying above a place in $S$. We shall assume that for the relevant places, we have {\it good reduction}. We do not need any subtle notion here and by this we simply mean that the reduced varieties $\bar W,\bar V$ are quasi-projective and of finite type over the residue field, and have the same dimension as  $W,V$. In turn, to define the reduced varieties, it suffices to take a model for $V,W$ by polynomial equations and reduce the equations modulo $v$. For large enough $S$ (depending on the models)  it is standard (and easy to prove)  that the dimensions are preserved.

Then  we define the set of $(k,S)$-integral points, or simply $S$-integral points,  $V(\O_{k,S})$ as the set of points $P\in V(k)$ such that,  for each place $v$ of $k$ not lying above a place in $S$, the reduction of $P\pmod v$ lies in the reduction of $V\pmod v$ (note that the reduction of $P$ at $v$  is  always well-defined as a point of $\P_n$  - because $v$ is a discrete valuation - with coordinates  over the residue field). Namely, we do not want that the reduction of $P$ lies in the complement of $V$ with respect to its projective closure.   In practice, e.g. when $V$ is given by an affine embedding, $V(\O_{k,S})$ consists of the points in $V(k)$ having coordinates   in $\O_{k,S}$.  As before, we may omit the reference to $k$ and speak of $S$-integral points meaning the points defined over $\overline\Q$ and being $(k',S)$-integral for every number field $k'$ over which they are defined.  We denote by $V(\O_S)$ the set of  these points.

Note that if   $V$ is projective then we have no condition and by $(k,S)$-integral points we recover the usual $k$-rational points.

We have:

\begin{thm}[Chevalley-Weil for integral points]
  \label{thm:CW_integral}
  Let $\pi: W \to V$ be a dominant morphism of 
  quasi-projective varieties defined over a number field $k$. Assume that there exists an embedding $\varphi: k \to \C$ such that the corresponding map $\pi_\C: W(\C)\to V(\C)$ is an unramified topological cover.
  Then, there exists a finite set $S$ of places of $\Q$ 
  such that if $P \in V(\O_{S})$ is an  $S$-integral point, $Q \in W(\bar k)$ and $\pi(Q) = P$ then  $Q$ is an $S$-integral point and the extension $k(Q)/k(P)$ is  ramified at most above places  of  $k(P)$ lying above  $S$. 
\end{thm}

We recover the former version from this one since, as noted above,  for a projective variety the $S$-integral points (resp. over $k$) are the usual $\overline\Q$-rational (resp. $k$-rational) points (for any $S$), since the complement of $V$ with respect to its closure is empty. 

From the fact that there exist only finitely many extensions of a given number field of fixed degree and unramified outside a fixed finite set of places (Hermite's theorem) we obtain the following, more classical, formulation:

\begin{cor}\label{cor}
Let $\pi:W\to V$ be a morphism of quasi projective varieties satisfying the hypotheses of the above theorem. Let $S$ be a finite set of places of $k$ containing the archimedean ones. Then there exists a finite extension $k'$ of $k$ and a finite set $S'$ of places of $k$, containing all those above $S$, with the following property: for every point $P\in V(\O_S)$ and every $Q\in W(\bar{k})$ with $\pi(Q)=P$, $Q$ lies in $W(\O_{S'})$.
\end{cor}
 
We note that this corollary, unlike the classical statements requiring $\pi$ to be an unramified cover (in the algebraic sense), is non trivial already for morphisms of degree one. Indeed, the induced continuous map $\pi_\C$ might be an homeomorphism  even when the morphism $\pi$ is not an isomorphism of algebraic varieties (see Example \ref{ex:cusp}).

\subsection*{Acknowledgements} We thank Julian Demeio and Ariyan Javanpeykar for discussions on conjugate varieties. AT was partially supported by Centro di Ricerca Matematica Ennio de Giorgi and is a memeber of GNSAGA of INdAM.

\section{Proof of the Chevalley-Weil Theorem(s).}\label{sec:proof}

We now offer the announced  
proof of the above versions of the Chevalley-Weil Theorem \ref{thm:CW} and Theorem \ref{thm:CW_integral}.  Actually we shall prove the theorems at the same time. (The proof of the former statement would result in a proof similar to that   presented in the book \cite{CZ} by the first and third authors).

We shall assume the embedding $\varphi:k\to\C$ is given and usually omit the reference to it  (and the corresponding subscripts) throughout. (See the next section for a proof that the topological assumption is independent of the choice of $\varphi$.) 

\begin{rem}  {\tt Effectivity}. It will be clear from the arguments that they are effective, in the sense that when an `effective presentation'  for the given varieties is given,  then one can determine  the relevant   finite sets of places  so that the conclusions hold.\footnote{This requires in particular effective versions  for the Nullstellensatz, and there are `good'  ones  in the literature. A somewhat different issue is to check effectively whether the topological assumptions hold. We believe this can be also done, though we do not have a reference, and we do not pause further on this aspect here.} 

\end{rem}

\medskip

\begin{proof}[Proof of Theorems \ref{thm:CW}  and  \ref{thm:CW_integral}]

As mentioned, it suffices to prove Theorem \ref{thm:CW_integral}, the former result being the special case when $V,W$ are projective.

\medskip

We start by observing that the cover $\pi: W\to V$ may be assumed to be Galois, i.e. connected and endowed with a  group of fiber-preserving automorphisms which is transitive on each fiber.  Indeed we prove the following

\begin{prop}
If Theorem  \ref{thm:CW_integral} holds for Galois covers, then it holds in general.
\end{prop}

\begin{proof}
Let   $\pi: W\to V$ be a dominant morphism satisfying the hypotheses of Theorem \ref{thm:CW_integral}. We want to construct a variety $W'$ and a dominant morphism $\pi': W'\to W$, possibly defined over a finite extension of the  number field $k$, such that the composition $\pi\circ \pi'$ is a Galois covering still satisfying the hypothesis of Theorem \ref{thm:CW_integral}. We then apply to that cover Theorem \ref{thm:CW_integral}, which is supposed to hold for Galois covers, and clearly we deduce the sought conclusion for the original morphism $\pi:W\to V$.

We want then to construct  the  Galois closure of the topological cover $\pi_{\C}: W(\C)\to V(\C)$. This is a topological space $\mathcal{W}$ endowed with a continuous map $ \mathcal{W} \to W(\C)$, which is a  cover of topological spaces,  such that the composed  map $\mathcal{W}\to W(\C)\to V(\C)$ is a Galois cover. 

\smallskip

In the first place, we note that the fibers of $\pi$ are finite and of the same cardinality, denoted  $d$, since $\pi_\C$ is a topological covering space.

Then, we let $W'$ be the a connected component of the algebraic variety consisting of the points $(x,y_1,\ldots,y_d)$, where $x\in V$ and $\{y_1,\ldots ,y_d\}$ is the fiber in $W$ above $x$, where $d=\deg\pi$ so the $y_i$ are distinct (recall that $\pi$ is a topological covering map).  Note that $W'$ is still quasi-projective and $W'(\C)$ is an unramified  cover of $V(\C)$ under the projection map $(x,y_1,\ldots,y_d)\mapsto x$, still denoted $\pi$. This $W'$ can be reducible but this will not matter in the present proof.  It is still a connected quasi-projective variety. 

\end{proof}

We now prove Theorem \ref{thm:CW_integral} assuming that the morphism $\pi: W\to V$ is a Galois cover. This means that 
we have a transitive action of a  group $G$  on the fibers of $\pi$ in $W$, without fixed points.  

For large enough  $S$ and any place $\nu$ of $\overline\Q$ not lying above $S$, we may reduce the varieties modulo $\nu$, and obtain  an algebraic cover $\bar\pi :\bar W\to\bar V$, Galois with Galois group $\bar G$, where $\bar\pi, \bar G$ are obtained as the reductions of the corresponding maps.  Here we do not need any subtle concept of reduction, as in the above preliminary  comments.

\medskip

Now we give an explicit and direct proof that,  at the cost of   enlarging suitably  the finite  set $S$, the reduction $\bar G$ still acts on $\bar W$ without fixed points. 

Indeed, we have already noted that since $\pi$ is a topological covering map $\pi:W(\C)\to V(\C)$, the group $G$ acts without fixed points (on the whole $W(\C)$).  Thus for $g\in G$, $g\neq 1$, the closed subvariety $W^g$ of $W$ defined by $g(x)=x$ is empty.

Now, we may write $W=W_1\setminus W_2$ where $W_1\supset W_2$ are closed algebraic varieties of a suitable projective space $\P_m$, both defined over a number field, still denoted  $k$ (and of course $W_2$ will be empty if we aim at Theorem  \ref{thm:CW}). 

Let $W_3$ be the union over $g\in G-\{1\}$, of the closures in $W_1$ of the subvarieties $W^g$ of $W$ defined by $g(x)=x$.   Our assumption means that, for each given $g\neq 1$,   the  equation $g(x)=x$  (which amounts to a certain system of polynomial equations)  has no   solution $x\in W(\C)$, hence $W_3(\C)\subset W_2(\C)$. 

Let  further $\J_i$, $ i=1,2$  be the (homogeneous) ideal defining $W_i$ in $\P_m$, and let $\I$ be the ideal generated by the (homogeneous) polynomials, denoted $f_1,\ldots ,f_q$, defining the union of  the conditions $g(x)=x$, $g\in G-\{1\}$, $x\in W$. 
 Note that the ideal  $\J_1+\I$ defines the variety $W_3$ (though it need not be the whole  ideal of that variety). 
 
 Since $W_3(\C)\subset W_2(\C)$, we have ({\it Nullstellensatz}) that $\J_2$ is contained in the radical of $\J_1+\I$. Hence if $f_1,\ldots ,f_q$ are  generators for $\I$, if $f_{q+1},\ldots ,f_r$ denote  generators for $\J_1$ and if $h_1,\ldots ,h_s$ denote generators for $\J_2$ (all homogeneous of positive degree in  $k[x_0,\ldots,x_m]$, so  for instance  the $h_i$ can be taken as the variables $x_i$ when $W_2$ is empty), there exists an integer $N>0$ such that we have identical relations
\begin{equation}\label{E.CW}
h_i^N=\sum_{j=1}^ra_{ij}f_j,\qquad i=0,\ldots ,s,
\end{equation}
for suitable homogeneous polynomials $a_{ij}\in k[x_0,\ldots ,x_m]$. Conversely, these equations imply that $W_3\subset W_2$.
(Note that when $V,W$  are  projective $W_2$ will be empty, in which case $\J_2$ would contain a power of the irrelevant ideal generated by the coordinates $x_i$.)

\medskip

%%%%%%%%%%%%%%%%%%%%%%%%%%%%%%%%%%%%%%%%%%%%%%%%%%%%%%%%%%%%%%%%%%%%%%%%%%%%%%%%%%%%%%%%%%%%%%%%%%%%%%%%%%%%%%%%%%%%%%%%%%%%%%%%%%%%%%%%%%%%%%%%%%%%%%%%%%%%%%%%%%%%%%%%%%%%%%%%%%%%%%%%%%%%%%%%%%%%%%%%%%%%%%%%%%%%%%%%%%%%%%%%%%%%%%%%%%%%%%%%%%%%%%%%%%%%%%%%%%%%%%%%%%%%%%%%%%%%%%%%%%%%%%%%%%%%%%%%%%%%%%%%%%%%%%%%%%%%%%%%%%%%%%%%%%%%%%%%%%%%%%%%%%%%%%%%%%%%%%%%%%%%%%%%
\begin{comment}Let $f_1,\ldots ,f_m \in k[x_0,\ldots,x_n]$ be homogeneous polynomials defining  $W_g$ in $\P_n$ (obtained by taking the equations defining $W$ and adding the equations corresponding to $g(x)=x$). 
%by the above we may assume that $f_j\in  k[x_0,\ldots ,x_n]$. 
The  Nullstellensatz allows to translate the statement $W_g(\bar k)=\emptyset$ into identical relations
\begin{equation}\label{E.CW}
x_i^N=\sum_{j=1}^ma_{ij}f_j,\qquad i=0,\ldots ,n,
\end{equation}
for  a suitable integer $N>0$ and suitable homogeneous polynomials $a_{ij}\in   k[x_0,\ldots ,x_n]$ depending on $g$. 
\end{comment}
%%%%%%%%%%%%%%%%%%%%%%%%%%%%%%%%%%%%%%%%%%%%%%%%%%%%%%%%%%%%%%%%%%%%%%%%%%%%%%%%%%%%%%%%%%%%%%%%%%%%%%%%%%%%%%%%%%%%%%%%%%%%%%%%%%%%%%%%%%%%%%%%%%%%%%%%%%%%%%%%%%%%%%%%%%%%%%%%%%%%%%%%%%%%%%%%%%%%%%%%%%%%%%%%%%%%%%%%%%%%%%%%%%%%%%%%%%%%%%%%%%%%%%%%%%%%%%%%%%%%%%%%%%%%%%%%%%%%%%%%%%%%%%%%%%%%%%%%%%%%%%%%%%%%%%%%%%%%%%%%%%%%%%%%%%%%%%%%%%%%%%%%%%%%%%%%%%%%%%%%%%%%%%%%%%%%%%%%%%%%%%%%%%%%%%%%%%%%%%%%%%%%%%%%%%%%%%

Since  the set  of all the coefficients of the polynomials $a_{i,j}$  is finite, there is a finite set  $\Sigma$  of places of $k$, such that for every $v\not\in \Sigma$ such coefficients are $v$-integers, and therefore the reductions modulo $v$ of the polynomials $f_1,\ldots, f_r$ generate  an ideal containing the reductions of the $h_i^N$.  We may assume that $S$ contains all the places of $\Q$ lying below $\Sigma$, so that this assertion holds for the places not above $S$.

Thus, we obtain that for all places $\nu$ of $\overline \Q$ not lying above a place in $S$ the reduction modulo $\nu$ of the variety $W_3$ is  contained in the reduction of $W_2$, and by construction we also have that  the reduction $\bar G$ of $G$ modulo $\nu$ continues to act without fixed points on the reduction $\bar W$ of the variety  $W$.

\medskip

Before going ahead we pause to note  that the assumptions   imply that $\pi$ is a {\it finite map}. Indeed, first of all the fact that it induces a topological covering space on the complex points implies of course that it has finite degree and in fact is quasi-finite (i.e., it has finite fibers). Then one can note that $\pi$ is also a proper morphism (for instance we note that  the topological assumption implies that the inverse imagine by $\pi$ of a compact set is compact). Then one may apply a well-known result (originally of Deligne)  asserting that these two properties imply that the map is finite. This   appears for instance in \cite[Corollary 12.89]{GW}, or \cite[Lemma 14.8]{H} at p.178, or also \cite[Lemma p.220]{D}.

\medskip

Let us  now take  $P\in V(\overline\Q)$ and $Q\in W(\overline\Q)$ in the preimage $\pi^{-1}(P)$, assuming that $P$ is $S$-integral. 

The next step is to prove that, if $S$ has been chosen large enough in advance, $Q$ is also $S$-integral. For this, we shall use that $\pi$ is finite.

By definition (see for instance \cite{Ha}, p. 84), $V$ can be covered by affine open subvarieties $X_1,\ldots ,X_a$ such that  for every  $i=1,\ldots ,a$,  $Y_i:=\pi^{-1}(X_i)$ is also affine and the ring  $k[Y_i]$ is integral over $k[X_i]$. 

 We let $X_i$ be defined by $\phi_i\neq 0$ in $V$, where $\phi_i$ are homogeneous polynomials in the projective coordinates of  a projective  space in which $V$ is embedded.

We may suppose that $S$ is so large that all such subvarieties, polynomials  and maps introduced so far have good reduction modulo every place not lying above $S$. We may also assume that  integral equations for generators of the $k[Y_i]$ over $k[X_i]$  have coefficients that are $S$-integers (indeed, we have only finitely many objects to take into account). 

\medskip

Let $\nu$ be a place of $\overline\Q$ not lying above any place of $S$.  Since the reduction of $P$ modulo $\nu$ lies in $V$, we have that $\phi_i(P)\not\equiv 0\pmod\nu$  for some $i\in\{1,\ldots ,a\}$, and in particular $\phi_i(P)\neq 0$ and $P$ lies in $X_i$, hence $Q\in Y_i$. The affine coordinates of $Q$ satisfy integral equation over the affine coordinates of $P$. Since the latter are $\nu$-integral by assumption, and since the relevant equations have $\nu$-integral coefficients, the affine coordinates of $Q$ have also $\nu$-integral coefficients. Then we may reduce modulo $\nu$ and obtain that the reduction of $Q$ lies in  the reduction of $W$ (i.e., not in the reduction of $W_2$). 
(Alternatively, we may use directly the `Valutative Criterion of Properness', i.e. \cite[Theorem 4.7]{Ha} p. 101.)  
Since this holds for all places $\nu$ in question, we deduce that $Q$ is an $S$-integral point of $W$.

\bigskip

To go toward the end of the proof,  let $v$ be a place of $k(P)$ such that the restriction of $v$ to $\Q$ is not in $S$; we want to prove that the field extension $k(P,Q)/k(P)=k(Q)/k(P)$ is unramified above $v$ (note also  that $k(P,Q)=k(Q)$ since $P=\pi(Q)$). Let $w$ be any  place of $k(P,Q)$ above $v$.

Let also  $\Gamma$ be the Galois group of the Galois closure $L/k(P)$ of $k(P,Q)/k(P)$. For $\gamma\in\Gamma$, we have that $\gamma(Q)\in W$ and $\pi(\gamma(Q))=\gamma(\pi(Q))=\gamma(P)=P$, so $\gamma(Q)$ is also in the fiber $\pi^{-1}(P)$. Hence, since $G$ acts transitively on the fibers of $\pi$, there exists $g=g_\gamma\in G$ such that $g(Q)=\gamma(Q)$. 

\medskip

Now, suppose by contradiction that $k(P,Q)/k(P)$ is ramified at $w$ and let $\Gamma'\subset \Gamma$ be  the inertia group over $k(P)$ of a place of $L$ above $w$. If we had $\Gamma'(Q)=\{Q\}$, namely if $Q$ lay in the fixed field of $\Gamma'$, then $k(P,Q)/k(P)$ would be in fact unramified at $w$. Therefore we can assume that there exists $\gamma\in \Gamma'$ such that $\gamma(Q)\neq Q$. Let $g=g_\gamma$ be as above, and note in particular that $g\neq 1$. On the one hand we have $\gamma(Q)\equiv Q\pmod w$, so $g(Q)\equiv Q\pmod w$.  

On the other hand, and here is a crucial point, $Q$ is $S$-integral, hence its reduction modulo $w$ lies in $W$. Hence, this would imply that $\bar G$, the the reduction of $G$ modulo $w$, has a fixed point in the reduction of $W$, namely the reduction of $Q$ modulo $w$, contradicting the former conclusion. 
\end{proof}

As already observed, Theorem  \ref{thm:CW_integral}  implies  Theorem \ref{thm:CW}.
For the latter,  in fact a  shorter argument applies. 
 
\medskip

\section{Unramified topological covers}\label{sec:topo}

We discuss here the condition that $\pi_\C$ is an unramified topological cover, comparing it with other related conditions, which are commonly assumed in other treatments of the theorem.

\bigskip

{\tt Comparison with \'etale morphisms}. (For this concept we refer e.g. to \cite[A.12.14]{BG} p. 580, or \cite[Section 5.4]{D} p. 232, or \cite{Ha}, p. 275.) We note that finite morphisms which are \'etale induce topological cover of the corresponding complex spaces; in the other direction:  

\centerline{{\it when $V$ is normal our topological condition implies that $\pi$ is \'etale. }}

\smallskip

To check this implication, note that the condition on $\pi_\C$ implies that all the fibers of $\pi$ have the same (finite) cardinality. Then, when $V$ is normal it follows, e.g. from  the paper \cite{LS} by Lang and Serre, that $\pi$ is algebraically unramified. Moreover in this setting, by \cite[Theorem 14.129]{GW} (essentially due to Chevalley), the map $\pi$ is flat and therefore (by definition) $\pi$ is \'etale.

However if $V$ is not assumed to be normal,  our  condition is generally weaker that that of being \'etale, 
as we show in the following (indeed very simple) Example.  

\begin{example}\label{ex:cusp}
Let $V$ be a plane cuspidal cubic, for instance given by the equation $y^2=x^3$,  and let $\pi: W \to V$ be its normalization (in this case $W=\A^1$ with the map $t\mapsto (t^3,t^2)$ to $V$), all defined over $\Q$. In this case, as it is easy to see,  the topological map $\pi_\C: W_\C \to V_\C$, via the unique embedding $\Q  \hookrightarrow \C$, is a homeomorphism (in fact a universal homeomorphism in the sense of \cite[\href{https://stacks.math.columbia.edu/tag/04DD}{Definition 04DD}]{stacks-project}). However $\pi$ is not \'etale, since the sheaf of relative differentials $\Omega^1_{W/V}$ is supported at the cusp; explicitly,  the pull-back of the ideal $(x,y)\Q[x,y]$ of the affine ring of $V$ (i.e., the ideal of the $(0,0)\in V$) is the ideal $t^2\Q [t]$ in the affine ring of $W$, so it is the {\it square} of the ideal of $0\in W$. 

Our version of the Chevalley-Weil Theorem (Corollary \ref{cor})  in this particular example contains a non-trivial, although very elementary, arithmetic consequence. It asserts for instance that if $(x,y)\in\Z^2$ is an integral point on that curve, then $x$ divides $y$ in $\Z$. 
\end{example}

\medskip

{\tt Normality}. We stress here that in general, the concept of ramification of a finite map does not behave nicely without normality conditions. 
Following the notation of the paper \cite{Hilb} by the first and third author, we say that a (quasi projective) algebraic variety is {\it algebraically simply connected} if every cover of its normalization without rational section ramifies somewhere. 

We mention, even if it is not strictly related to our situation, that a normal complex variety is algebraically simply connected if and only if its topological fundamental group has no subgroups of finite index (see for instance  \cite[Proposition 1.1]{Hilb}). However, for non-normal varieties both directions fail: in one case the nodal cubic is not topologically simply connected, having fundamental group isomorphic to $\Z$, but its normalization is $\P_1$. In the other direction, examples constructed by Catanese show that there exist quotients of products of curves that are simply connected but their normalization has fundamental group $\Z^2$. (See footnote 10 in \cite{Hilb}: one of Catanese's  examples starts from a product $\P_1\times E$ where $E$ is an elliptic curve, taking the quotient by
the equivalence relation $(t_1, p_1) \sim (t_2, p_2)$ if and only if $ t_1 = t_2 = 0$ and $p_2 = -p_1$ or $(t_1, p_1) = (t_2, p_2)$.
We obtain a variety $X$ which may be shown to be simply connected. Its normalization is just the original
$\P_1\times  E$, whose fundamental group is $\Z^2$.) \medskip

\medskip

{\tt Open maps}. A key property of topologically unramified covers of complex algebraic varieties is that they are open maps: this means that open sets in the complex topology have open images. 
In the other direction,  any proper open map between complex connected manifolds whose fibers have the same cardinality is a  topological cover. 
However not all maps between complex varieties are open, even assuming that the cardinality of the fibers are constant. Here is an example:

\begin{example}\label{ex:node}
  Let $V$ be a plane nodal cubic (for instance defined by $y^2=x^3+x^2$), let $V'$ be its normalization (in this case $V=\A^1$ with map $t\mapsto (t^2-1,t^3-t)$) and let $W = V' \setminus \{ P \}$ where $P$ is one of the two preimages of the node (e.g., $W=\A^1-\{1\}$). Then $\pi: W \to V$ is an injective map which is not open.
\end{example}

In the above example, the morphism was not a finite one. If one takes  $W=V'$ (without removing the point $P$), the corresponding map is finite and non open, but it is not true that the fibers have the same cardinality. 

 Note that in the topological setting (i.e. outside the realm of algebraic varieties), a map can be continuous, surjective, proper with  fibers of constant cardinality without being open, in particular without being a topological cover. Here is an example:
 
 \begin{example}\label{ex:disk}
 Take for $X$ the disjoint union of the open unit disk $\Delta$ with a single point $P$, where $\{P\}$ is open; set $Y=\Delta$ and the map $X\to Y$ sends $z\mapsto z^2$ for $z\in\Delta$ while $P\mapsto 0$. It is proper, surjective and each of its fibers contains exactly two points. However it is not open, since the image of the open set $\{P\}$ is not an open set. 
 \end{example}

We now show, however,  that this phenomenon never occurs for algebraic varieties:

\begin{prop}\label{prop.criterion}
Let $V,W$ be complex irreducible quasi projective algebraic varieties and let  $\pi:W\to V$ be a finite morphism whose fibers have constant cardinality. Then $\pi:W\to V$ is a topological cover. In particular, it is an open map.
\end{prop} 

Note that, in view of Example \ref{ex:cusp}, it is not necessarily true that $\pi$ is unramified.  

\begin{proof}
Let $P\in W$ be a (complex) point; we shall construct local analytic sections over $P$.
The question is local, so we can consider a Zariski neighborhood of $P$ and  suppose that $V,W$ are affine. Also, after refining such a neighborhood, we can suppose that $W$ is defined in $V\times \A^1$ by a polynomial equation
\begin{equation*}
y^d+a_1(x)y^{d-1}+\ldots+a_d(x)=0
\end{equation*}
where $a_1,\ldots, a_d$ are regular functions on $V$, and $d=\deg \pi$ (indeed, it suffices to add to the ring $\C[W]$ one more function separating the fibers of $P$). Also, the discriminant of the above polynomial never vanishes. We can view the $d$-uple $(a_1,\ldots,a_d)$ as a regular map $W\to \A^d$, whose image lies in the complement of the zero set of the discriminant of the generic polynomial $y^d+t_1y^{d-1}+\ldots+t_d$. 
Denote by $S$ this open set in $\A^d=\C^d$ and by $p_0=(a_1(P),\ldots,a_d(P))\in S$ the image of $P$. By the implicit function theorem, there exist a neighborhood $S_0$ of $p_0$ in the complex topology and local analytic functions $f_1,\ldots,f_d:S_0\to \C$ such that for all ${\underline{t}}=(t_1,\ldots,t_d)\in S_0$ and every $i=1,\ldots,d$,
\begin{equation*}
f_i^d(\underline{t}) +t_1 f_i^{d-1}(\underline{t})+\ldots +t_d=0.
\end{equation*}
By continuity of the map $(a_1,\ldots,a_d):W\to \C^d$, the pre-image under such a map of the neighborhood  $S_0$ of $p_0$ is a neighborhood of $P$ in $W$. On such a neighborhood we can define the $d$ sections $V\ni x\mapsto (x,f_i(a_1(x),\ldots,a_d(x)))\in W$. This concludes the proof.
\end{proof}

{\tt Change of the field embedding}.
In the statement of Theorem \ref{thm:CW} we assume that the base change of the morphism  $\pi$,  through the embedding $\phi$ of $k$ in $\C$,  is a topological  cover. Given that the fibers of $\pi$ have constant cardinality, this is equivalent to requiring that $\pi_\C$ is an open map in the $\C$-topology. {\it A priori}  this condition depends on the choice of an embedding of the field of definition into $\C$. We give here a sketch of the proof that this is not the case, i.e. being an unramified topological cover can be checked with respect to \emph{any} embedding of the field of definition. Maybe this is known
but we are  including a statement and proof  here for a lack of reference. 

If $X$ is an algebraic  variety defined over a number field $k$, explicitly embedded in some projective space by equations  defined over $k$, and if $\sigma\in G_\Q:={\rm Gal}(\overline\Q/\Q)$, we speak of a {\it conjugate variety $X^\sigma$} meaning the object defined by the equations obtained by applying $\sigma$ to the coefficients of the equations defining $X$.  The conjugate variety is defined over the image $\sigma(k)$ of the number field $k$.

It  will make no difference for our results to enlarge the number field $k$,  and in particular   to assume that $k/\Q$ is a Galois extension; so we adopt tacitly this convention from now on.

 Given one embedding $\phi:k\hookrightarrow \C$, all the other embedding are of the form $\phi^\sigma:=\phi\circ \sigma$, where $\sigma\in {\rm Gal}(\overline\Q/\Q)$. Extending $\sigma$ to an automorphism of $\C$, we obtain a bijection between the complex points of $X$ and of $X^\sigma$. Unless this extension is the complex conjugation or the identity, such a bijection will not result in a continuous map.

Now, given a morphism $\pi:Y\to X$ of algebraic varieties over the number field $k$, an embedding $\phi:k\hookrightarrow \C$ and a Galois automorphism $\sigma$, we obtain the following commutative diagram
\begin{equation}\label{diag}
    \xymatrix{Y \ar[r]^\sigma \ar[d]_\pi & Y^\sigma \ar[d]^{\pi^\sigma} \\
              X \ar[r]^\sigma & X^\sigma}
    \end{equation}
where the horizontal arrows are in general discontinuous maps. A priori, it is not obvious that whenever the vertical arrow $\pi$ is a topological cover, so is $\pi^\sigma$. This is the content of the following 

\begin{prop}\label{prop:conj}
Let $\pi_k: Y \to X$ be a finite dominant morphism between quasi-projective varieties defined over a number field $k$. Fix an embedding $k\hookrightarrow\C$ and let $\pi: Y(\C) \to X(\C)$ be the corresponding map between complex points. For an automorphism $\sigma\in {\rm Gal}(\overline\Q/\Q)$ let $\pi^\sigma: Y^\sigma (\C) \to  X^\sigma(\C)$ be the corresponding map between the complex points  as above. Then $\pi$ is an unramified topological cover of the complex points  if and only if $\pi^\sigma$ is an unramified cover.
\end{prop}

\begin{proof}
The condition of being a finite map is preserved by extension of scalars; so if $\pi_k$ is a finite morphism of algebraic varieties over $k$, which we are assuming, also the morphisms $\pi_\phi$ of complex algebraic varieties are finite, for every chosen embedding $\phi: k\to \C$. So in particular $\pi,\pi^\sigma$ are finite maps. Now, if one of them is a topological cover, then its fibers have the same cardinality. Clearly, this fact must hold also for its conjugate. By Proposition \ref{prop.criterion} above, the latter map is a topological cover.
\end{proof}

Looking back at the diagram \eqref{diag}, we have already remarked that the horizontal bijective arrows are not continuous, so in general do not send the  neighborhoods of a point to neighborhoods of its conjugate. However, we prove the following weaker fact, which might be of independent interest. 
 
\begin{prop}
  \label{prop:open}
Let $X$ be an algebraic  variety defined over a number field $k$.  Let $P$ be a smooth algebraic point of $X$, let $\sigma \in \Gal (k(P)/\Q)$ and $Q = P^\sigma$ be the conjugate in the corresponding conjugate variety $X^\sigma$. Then, for any open set $V_P$ containing $P$ and $V_Q$ containing $Q$ there exist an extension of $\sigma$ to $\overline\Q$, an open neighbourhood $U_P\subset V_P(\C)$ of $P$ in $X$, in the complex topology,  that contains a dense set $\{ P_i \}_{i \in \N}$ of algebraic points of $X$ and an open neighbourhood $U_Q \subset V_Q$ of $Q$ in $X^\sigma(\C)$, such that $P_i^\sigma \in U_Q$ for all $i \in I$ and the set $\{P_i^\sigma: i\in \N\}$ is dense in $U_Q$.
\end{prop}

\begin{proof}
  The question is local so we may assume that $X$ is affine. By Noether normalization \cite[\href{https://stacks.math.columbia.edu/tag/0CBH}{Lemma 0CBH}]{stacks-project} there exists a finite morphism $f: X \to \A^d$ where $d = \dim X$. Locally the map $f$ is given by (a general) projection, so we can assume the existence of an absolutely irreducible polynomial $F(y,x_1,\dots,x_d)$, monic in $y$ of degree $m=\deg f$, and  defined over $k$,  such that in a neighbourhood of $P$ the variety $X$ is given by $F = 0$ and the map $f$ is the projection to $(x_1,\dots,x_d)$. Also, we may assume that the discriminant of $F$ with respect to $y$ does not vanish at $P$, since $P$ is a smooth point of $X$. 
  
  The assumption on the discriminant  entails that the projection map is a topological cover in  $f^{-1}(f(U))$ for $U$ a small open neighbourhood of $P$ (not  projecting to points in the zero locus of the discriminant). 
  Hence if $D$ is a small enough disk around $f(P)$ in the affine space $\A^d$, its inverse image by $f$ will be a disjoint union of neighbourhoods of the  points $(z_\mu, f(P))$, where $z_\mu$ ranges over the $m$ distinct roots of $F(y,f(P))=0$. We choose then a small enough $D$ and we  take as $U_P$ the corresponding neighbourhood of $P$. A suitable choice of $D$ will ensure that $U_P$ is contained in the given open set $V_P$.

  Then $D$ contains a dense set of points of the form $f(P) + \underline{t}$ where $\underline{t} \in \A^d(\Q)$,  $| \underline{t} |$ is small, and     
   the preimage of $D$ in $X$ contains an open neighbourhood $U_P$ of $P$   containing  a dense set of  algebraic points $P_i$ that verify $f(P_i) = f(P) + \underline{t_i}$, $i\in\N$.   
  
  We can in fact refine this construction. Suppose the sequence of $t_i$ has been chosen as above. Then we deform it inductively  as follows. Suppose that  $t'_1,\ldots ,t'_{n-1}$ have been determined, and put $P'_i=(y_i,f(P)+t'_i)\in X(\C)$, for suitable algebraic numbers $y_i$. Then we  choose  $t'_n$ very near to $t_n$ (say $|t_n-t'_n|<1/n$) and such that 
  the polynomial $F(y,f(P)+t'_n)$ is irreducible over the Galois closure over $\Q$ of the field $k(P,y_1,\ldots ,y_{n-1})$. This is possible by Hilbert Irreducibility Theorem (in a refined form easy to prove, asserting  density of good specialisations, there are many such ones, see e.g. \cite{BG}). 
  
  Now we choose $y_n$ as the root of $F(y,f(P)+t'_n)$ nearest to the coordinate $y(P)$; actually, if $D$ is small (as prescribed) we shall have $P_n:=(y_n,t'_n)\in U_P$.  This completes the induction step. 
  
  Plainly the new sequence continues to fulfil the former properties, so we use it in place of the former, omitting the dash in the notation.

  \medskip

 To go ahead, by conjugating with $\sigma$ we obtain the following diagram:
  \[
    \xymatrix{X \ar[r]^\sigma \ar[d]_f & X^\sigma \ar[d]^{f^\sigma} \\
              \A^d \ar[r]^\sigma & \A^d}
    \]
   
   Note that $f^\sigma$ will still be a  finite map and a covering map above a neighbourhood of $f(Q)=f(P)^\sigma$. 
   
     By induction on $n$, we can choose an extension of $\sigma$ to the field $k(P,y_1,\ldots ,y_n)$   so that the conjugates $P_i^\sigma$  are dense in a small neighbourhood of $P^\sigma=Q$. In fact, let us again suppose we have extended $\sigma$ to $k(P,y_1,\ldots, y_{n-1})$ in a certain way.  Now, in the first place we have automatically that $f(P_n)^\sigma=(f(P)+t_n)^\sigma=f(Q)+t_n$, since the $t_n$ are rationals and so $t_n^\sigma=t_n$.  Therefore, since $F^\sigma(y_n^\sigma,f(P_n)^\sigma)=(F(P_n))^\sigma=0$, we deduce that $y_n^\sigma$ will be  near some zero of $F^\sigma(y,f(Q))$. These zeros include $y(Q)=y(P)^\sigma$, so one of the zeros of $F^\sigma(y,f(P_n)^\sigma)$ will be near $y(Q)$. Since the last polynomial is by construction irreducible over $k(P,y_1,\ldots ,y_{n-1})$, we can compose $\sigma$ on the left with an automorphism fixing $k(P,y_1,\ldots ,y_{n-1})$ and sending $y_n$ to the root  of $F^\sigma(y,(f(P_n)^\sigma)$    nearest  to $y(Q)$. We replace the previously chosen $\sigma$ with this new one (and the action will be the same on the field generated by the first $n-1$ points). For $n\to\infty$ then $y_n^\sigma$ will lie in a prescribed neighbourhood of  $y(Q)$ and be dense there, since the $f(P_n)^\sigma$ are dense in a disk around $f(Q)$. Up to further shrinking the disk $D$ (and hence $U_P$) we can always choose the points to all lie in the given $V_Q$.
   
  This proves the assertion of the proposition.
  \end{proof}

\begin{rem}
  We note that the assumption that $P$ is algebraic can be relaxed as follows: one can consider the field generated by $P$ in $k$ and then use the same argument to produce points $P_i$ in a neighbourhood of $P$ that are algebraic over this new field.
\end{rem}

\section{An application}\label{sec:appl}
As an application of Theorem \ref{thm:CW_integral} we show how to recover the following result of Darmon and Granville on generalized Fermat equations. The original argument uses a Theorem of Beckmann \cite{Beck}, while here we follow the notes \cite{C}.

\begin{thm}[{\cite[Theorem 2]{DG}}]
  \label{th:DG}
Let $(p,q,r)$ be a hyperbolic triple of positive integers, i.e. such that $1/p + 1/q + 1/r < 1$. Then for every non-negative integers $a,b,c$ there exist only finitely many solutions $(x,y,z) \in \Z^3$ with $\gcd(x,y,z) = 1$ to the generalized Fermat equation
\begin{equation}\label{eq:Fermat}
  a x^p + b y ^ q = c z^r.
\end{equation}
\end{thm}

We first note that integral solutions to equation \eqref{eq:Fermat} correspond to integral points of the affine surface $\S'$ defined in $\A^3$ by the same equation. The coprimality condition, $\gcd(x,y,z) = 1$, is equivalent to the fact that we are only considering points $(x,y,z)$ in the quasi-projective variety $\S = \S' \setminus \{\underline{0}\}$ that are integral with respect to the origin $\{ \underline{0} \}$. In fact given a point $P \in \S$ the fact that its reduction modulo a prime $\nu$ is not $\underline{0}$ is equivalent to the fact that $\nu$ does not divide the $\gcd$ of the affine coordinates of $P$. We also note that in this situation it is crucial that our previous statements apply to quasi-projective varieties that are not necessarily affine like the surface $\S$.

Following Darmon and Granville, one considers the morphism $\beta: \S \to \P_1$ given by
\[
  \beta(x,y,z) := \dfrac{ax^p}{cz^r}.
\]

The map $\beta$ has multiple fibers over $\{ 0, 1, \infty \}$ (we will discuss this more in detail in the proof of Proposition \ref{prop:DG}) but one can remove these multiplicities via a, necessarily ramified, cover of the base $\P_1$. More precisely, there exists a ramified cover $\pi: \calC \to \P_1$ ramified precisely over $\{ 0 , 1, \infty \}$ and with ramification orders $(p,q,r)$ (e.g. one can construct such cover from non-trivial finite quotients of the triangle group $T(p,q,r)$, see \cite[Section 5]{C}). Then, letting $\calF$ be the \emph{normalization} of the fiber product $\S \times_{\P_1} \calC$ of $\pi: \calC \to \P_1$ and $\beta: \S \to \P_1$, we can define a finite dominant map $\bar\pi: \calF \to \S$. Then the map $\bar\pi$ satisfies the hypotheses of Theorem \ref{thm:CW_integral}. This is the content of the following proposition.

\begin{prop}[{see \cite[Lemma 5.10]{C}}]
  \label{prop:DG}
  In the above setting, the map $\bar\pi_\C: \calF_\C \to \S_\C$ is an unramified topological cover.
\end{prop}

Note that we could not apply Theorem \ref{thm:CW_integral} to the natural map $\S \times_{\P_1} \calC \to \S$, since the latter is ramified over $\beta^{-1}(\{0,1,\infty\})$. This is the reason why in \cite{DG} the authors have to apply a different argument, using \cite{Beck}.

We now show how an application of the Chevalley-Weil Theorem \ref{thm:CW_integral} implies Theorem \ref{th:DG}. 
\begin{proof}[Proof of Theorem \ref{th:DG}]
  We have already shown how coprime solutions of equation \eqref{eq:Fermat} correspond to integral points of $\S$. By Proposition \ref{prop:DG} we can apply Theorem \ref{thm:CW_integral} to $\bar\pi: \calF \to \S$ and obtain that the integral points of $\S$ lift to integral points of $\calF$ defined over a fixed ring of $S$-integers of a number field $k$. The image of an $S$-integral point of $\calF$ in the curve $\calC$ is by construction a rational point of $\calC$ defined over $k$. 

   Now we note that, since the integers $p,q$ and $r$ satisfy the condition $1/p + 1/q + 1/r < 1$, the Riemann-Hurwitz formula shows that the genus of $\calC$ is at least 2. Hence, the celebrated Faltings' Theorem implies that the $k$-rational points of $\calC$ are finite. Therefore, the rational points, and hence the integral points, of $\calF$ are contained in a finite number of fibers of the map $\calF \to \calC$. Then, the same holds true for the integral points of $\S$, namely they are contained in finitely many fibers of the map $\beta: \S \to \P_1$. But the fibers of $\beta$ are all isomorphic to the multiplicative group $\G_m$ and therefore can contain only finitely many points with integral coprime coordinates, concluding the proof.
\end{proof}

We end this section with a proof of Proposition \ref{prop:DG}. Before giving the detailed proof, we note that Proposition \ref{prop:DG} is a natural extension to dimension two of the following fact on covers of curves; it is one formulation of the so-called Abhyankar's lemma (see \cite{SeSem} for a discussion in the context of covers of the plane):\smallskip

{\it Let $\mathcal{X}$ be a (complex) projective algebraic curve and $f\in\C(\mathcal{X})$ be a non-zero rational function whose divisor is divisible by an integer $n>1$. Let $\mathcal{Y}$ be a normal projective  model of the function field $\C(\mathcal{X})(\sqrt[n]{f})$. Then the natural morphism $\mathcal{Y}\to \mathcal{X}$ is unramified}.

\smallskip

Note that the equation $y^n=f(x)$ defines in $\P_1\times \mathcal{X}$ a {\it singular} model of the same function field of $\mathcal{Y}$. The same happens in our two-dimensional situation.

\begin{proof}[Proof of Proposition \ref{prop:DG}]
 The question is local, so we argue locally (in the complex topology) and we will omit the subscripts, assuming that in this proof we only refer to the complex manifolds associated to the algebraic varieties considered above. Let $\bar \beta: \calF \to \calC$ be the projection to $\calC$ so that, by definition, $\pi \circ \bar \beta = \beta \circ \bar \pi$. This is made explicit in the following diagram.
 \[
   \xymatrix{ \calF \ar[r]^{\bar{\pi}} \ar[d]_{\bar\beta} & \S \ar[d]^{\beta} \\ \calC \ar[r]^{\pi} & \P_1}
 \]
 
 We already observed that the fibers of $\beta$ are all isomorphic to $\G_m$; in fact $\beta$ restricted to $\P_1 \setminus \{ 0,1, \infty \}$ is a (principle) $\G_m$-bundle (but it is not over the whole $\P_1$). In particular, given $s \in \S$ such that $\beta(s) \notin \{ 0,1, \infty\}$, there exists a neighbourhood $U$ of $\beta(s)$ such that $\pi: \pi^{-1}(U) \to U$ is a topological cover. In this case, the surface $\calF$ is locally defined as $\beta^{-1}(U) \times_U \pi^{-1}(U)$. Since the latter is a smooth complex space it follows that $\bar \pi: (\pi \circ \bar\beta)^{-1}(U) \to \beta^{-1}(U)$ is a topological cover (the pull-back via $\beta$ of $\pi$ restricted to $\pi^{-1}(U)$).\medskip

 On the other hand, when $\beta(s) \in \{0,1,\infty\}$ the fiber-product of $\S \times_{\P_1} \calC$ is singular above $\beta(s)$ (but its normalization $\calF$ is smooth). We will show that also in this case, the map $\bar \pi$ is a local biholomorphism by making use of the interplay between the multiplicities of the fibers of $\beta$ over $\{ 0, 1, \infty \}$ and the ramification of the map $\pi$.

 We will give details for the case when $\beta(s) = 0$, since the other two cases are dealt with in a similar way. In this case we let $s \in \S$ such that $\beta(s)= 0$. We can consider as local parameters around $s$, in the surface $\S$ embedded in $\A^3$, the regular functions $x$ and $z - z(s)$. In fact the function $y$ can be expressed as
 \begin{equation}\label{eq:y}
  y =  \left( \dfrac{c z^r - a x^p}{b} \right)^{1/q}.
 \end{equation}
 Here the $q$-th root function is well-defined locally at $s$ and we assume that a choice of a branch compatible with $y(s)$ has been made.

 Now we consider a point $f \in \calF$ that lies in the fiber above $s$, i.e. such that $\bar\pi(f) = s$, and let $\gamma = \bar\beta(f)$ be the corresponding point in the curve $\calC$. Since the (Galois) cover $\calC \to \P_1$ is ramified at 0 of order $p$, there exists a local parameter $t$ at $\gamma$ in $\calC$ such that $t^p = \bar\pi^*\beta$. In particular the cover $\calC \to \P_1$ will be given locally by $t \mapsto t^p$. We stress here that the relation $t^p = \bar\pi^* \beta$ defines a singular (in fact not normal) variety. On the other hand the local ring of $\calF$ over the local ring of $s$ in $\S$ is integrally closed and generated by the element $t$. In other words the surface $\calF$ is defined (birationally) by adding the function $t$ to the local parameters $x$ and $z-z(s)$. Note that, by definition, the function $t$ satisfies the relation
 \begin{equation}\label{eq:t}
  t = \dfrac{x}{\left(cz^r/a\right)^{1/p}}.
 \end{equation}
 As before the $p$-th root function is well defined locally and we assume a choice of its branch has been made accordingly to the choice of the point $f \in \calF$ lying over $s$. We stress that in this situation there are $p$ possible choices for the point $f \in \calF$ such that $\bar\pi(f) = s$ and $\bar\beta(f)=\gamma$. Such choices correspond to the branches of the $p$-th root of the denominator of $t$ in the equation \eqref{eq:t}. On the other hand, if we would have considered the surface $S \times_{\P_1} \calC$ we would have had only one choice. 

 Finally the functions $t$ and $w:= \bar\pi^*(z - z(s))$ are local parameters at $f$ in $\calF$, in particular we can explicitly compute the following expressions
 \[
   z = w + z(s) \qquad x = t\left( \dfrac{c z^r}{a} \right)^{1/p},
 \]
 and similarly for $y$ using equation \eqref{eq:y}. This yields that the function $\bar\pi$ can be defined as $(t,w) \mapsto (x,y,z)$ and therefore it is a biholomorphism. In particular this shows that $\bar\pi$ is unramified as wanted.
\end{proof}

\end{document}